\magnification=\magstep1
\baselineskip=16pt
\font\bgrm=cmbx12
\vglue1cm

\centerline{\bgrm Semigroups in finite von Neumann algebras}

\bigskip

\centerline{\bf Gilles Cassier}

\bigskip

\centerline{\it  Dedicated to the memory of Professor B\'ela Sz\H okefalvi-Nagy}

\bigskip

{\bf Abstract.} Let $M$ be a finite von Neumann algebra. In the first
part, we give asymptotic results about $M$-stable sequences of
weak*-continuous mappings which are related with operators belonging to $M$. In
the second part, we extend, by a shorter way, similarity results given in
[CaFa2] to unbounded semigroups of operators contained in a finite von
Neumann algebra.

\bigskip

{\bf I. Introduction and preliminaries}

\medskip

Let $H$ be a separable complex Hilbert space and let $B(H)$ be the algebra of
bounded linear operators acting on $H$. The ultra-weak topology of $B(H)$ is the weak*
topology (in the sequel we will shorten weak* to w*) that comes from the
well known duality $B(H)=(C_{1}(H))^{\ast }$, where $C_{1}(H)$ is the Banach
space of trace class operators on $H$ endowed with the trace norm (see
[Dix]). A von Neumann algebra acting on $H$ is by definition an
ultra-weakly closed *-subalgebra of $B(H)$. Such a von Neumann algebra $M$
is finite if it admits a faithful normal trace $\tau $, which means
that $\tau $ is an ultra-weakly continuous linear functional on $M$ satisfying:

1) $\tau (AB)=\tau (BA)$ for any $A,B\in M$;

2) for any positive element $A$ in $M$, we have $\tau (A)\geq 0$ and $\tau
(A)=0\Longrightarrow A=0$.

We denote by ${\cal T}(M)$ the set of all faithful normal traces acting
on $M$. A good example of a finite von Neumann algebra is the w*-algebra
generated by the left regular representation of a countable discrete group.
We will denote by $M_{\ast }$ the predual of $M$.
For any subset ${\cal F}$ of $M$, we shall denote by ${\cal F}^{\prime }$
the family of operators commuting with every element of ${\cal F}$.

Let $B(M)$ denote the algebra of bounded linear operators acting on $M$, and
let $B_{w}(M)$ stand for the algebra of operators $T\in B(M)$ wich are weak*-continuous. 
Recall that $\phi \in B_{w}(M)$ if and only if $\phi $ is the
adjoint of a bounded linear operator acting on the Banach space $M_{\ast }$ (see
for instance [BCP]). For any $\phi \in B_{w}(M)$, let $\phi_{\ast }$ denote
the uniquely determined operator whose (Banach space) adjoint is $\phi $, that is $(\phi_{\ast })^{\ast }=\phi $. 
For more details on von Neumann algebras, we
refer the reader to [Dix] and [Sak].

As usual $[A_{i,j}]_{1\leq i,j\leq n}\in {\cal M}_{n}(B(H))$ denotes the $n\times n$ matrix which
 acts on the orthogonal sum of $n$ copies of $H$; its entries are operators acting on $H$. We remind the reader that 
${\cal M}_{n}(M)$ inherits a unique structure of von Neumann algebra. Let 
$\psi $ be a linear mapping from $M$ into itself, we define $\psi_{n}:{\cal M}_{n}(M)\to {\cal M}_{n}(M)$ 
by $\psi _{n}([A_{i,j}]_{1\leq i,j\leq n})=[\psi(A_{i,j})]_{1\leq i,j\leq n}$. 
 We call $\psi $ $n$-positive if $\psi_{n}$ is positive (that is positive operators are transformed into positive ones)
 and we call $\psi $ completely positive
if $\psi $ is $n$-positive for all $n$.

We proved in [CaFa2] that a power bounded operator $T$ in a finite von Neumann
algebra $M$ is similar to a unitary element of $M$ if and only if $T^{n}x\not\to0$ 
for any $x\in H\backslash\{0\}$ ($T$ is said to be of class $C_{1\cdot}$ in the Sz.-Nagy--Foias terminology). 
We will extend this result into two
directions. On the one hand, we will consider general semigroups. On the
other hand, we will work with operators which are not necessarily power
bounded. To achieve this, we have to find a proper framework, which will allow
short and well adapted methods.

In similarity problems, the idea of using limits in the sense of Banach  comes
from B. Sz.-Nagy [Nag]. In the sequel, we frequently use this idea.
Recall that
a Banach limit is a state, that is a linear functional $L$ with $\| L\| =L({\bf 1})=1$,
acting on the classical space $\ell^{\infty }$ of all complex bounded sequences
and satisfying $L((u_{n+1}))=L((u_{n}))$. A bounded sequence $(u_{n})_{n\geq1}$ 
is said to be almost convergent to a complex number $c$ if
$$\lim_{n\to\infty }\sup_{k\in {\bf N}}\left| {1\over n+1}
\sum_{i=k}^{k+n}u_{i}-c\right| =0.$$
Lorentz proved in [Lor] that $(u_{n})_{n\geq 1}$ is almost convergent to $c$
if and only if for every Banach limit $L$ we have $c=L((u_{n}))$. 
A sequence $(u_{n})_{n\geq 1}$ is said to be strongly almost convergent to $c$ if the
sequence $(\left| u_{n}-c\right| )_{n\geq 1}$ is almost convergent to $0$.
We will say that a sequence $(\phi _{n})_{n\geq 1}$ of operators in $B_{w}(M) $
 is weakly almost convergent to $\phi \in B_{w}(M)$ if
 $[l,\phi_{n}(X)-\phi (X)]$ almost converges to $0$ for any $(l,X)\in M_{\ast }\times M$.

\medskip

{\bf Definition 1.1.} A mapping $p\colon {\bf N}\to (0,\infty)$
is called a {\it gauge} if there exists $c_{p}>0$ such that
the sequence $p(n+1)/p(n)$ is strongly almost convergent to $c_{p}$.
 Moreover, if in addition the sequence $c_{p}^{n}/p(n)$
strongly almost converges to $1$, then we say that $p$ is a
{\it regular gauge}.

\medskip

We will say that a sequence $(T_{n})_{n\geq 1}$ of operators, acting on a Banach space,
is {\it dominated} by a gauge $p$ if $\left\| T_{n}\right\| \leq p(n)$ holds for
every positive integer $n$.
 We follow [Ker] in saying that $(T_{n})_{n\geq1}$ is
{\it compatible} with a gauge $p$ if in addition the sequence $\left\|T_{n}\right\| /p(n)$
 does not almost converge to $0$.
An operator 
$T$ is dominated by (compatible with)  $p$ if the sequence $(T^{n})_{n\geq 1}$ is dominated by (resp. compatible with)
 $p$. Finally, a family ${\cal F}$\ of operators is called dominated by (compatible with)
 $p$ if each operator in ${\cal F}$\ is dominated by (resp. compatible with)  $p$. For some
recent contributions in this area, we refer the reader to [Ker], [Ker1],
[Ker2], [Ker3], [Ker4], [Ker5] and [KeM\"{u}].

Assume that $p$ is a gauge and $T\in M$ is dominated by $p$. Given a Banach limit $L$, let us introduce
 the (bounded, linear) operator $E_{L,T}$, acting on $M$, by setting
$$[l,E_{L,T}(X)]=L(\{[l,T^{\ast n}XT^{n}]p(n)^{-2}\}_{n\geq 1})$$
for any $(l,X)\in M_{\ast }\times M$. The following proposition summarizes
some useful properties of the operator $E_{L,T}$.

\medskip

{\bf Proposition 1.2.} {\it Let $T$ be an element in a von
Neumann algebra $M$ acting on a separable Hilbert space $H$.
 Assume that $T$ is dominated by a gauge $p$.
Then, for any Banach limit $L$, we have}

(i) $E_{L,T}$ {\it is a completely positive mapping};

(ii) $E_{L,T}(T^{\ast }XT)=c_{p}^{2}E_{L,T}(X)$ {\it for any }$X\in M$;

(iii) {\it if $A,B\in M$ commute with $T$, then we
have $E_{L,T}(A^{\ast }XB)=A^{\ast }E_{L,T}(X)B$ for any $X\in M$};

(iv) $T^{\ast }E_{L,T}(X)T=c_{p}^{2}E_{L,T}(X)$ {\it for any } $X\in M$;

(v) {\it there exists $\rho _{L}(p)\in [0,1]$ such that $E_{L,T}\circ E_{L,T}=\rho _{L}(p)E_{L,T}$};

(vi) {\it moreover, if $M$ is a finite von Neumann algebra,
then the mapping $E_{L,T}$ belongs to $B_{w}(M)$}.

\medskip

{\bf Remark 1.3.}  If $T\in M$ is compatible with a gauge $p$, then the
spectral radius $r(T)$ satisfies $r(T)=c_{p}$ (see [Ker)]).

\medskip

{\bf Proof.}   (i) Let $[X_{i,j}]_{1\leq i,j\leq p}$ be a positive $p\times p$
matrix whose entries are operators in $M$ and let $x_{1},...,x_{n}$ be
vectors in $H$. For any $(i,j)\in \{1,..,p\}^{2}$, we define the linear
functional $l_{i,j}$ acting on $M$ by setting $l_{i,j}(X)=\langle Xx_{j},x_{i}\rangle$.
It is obvious that $l_{i,j}\in M_{\ast }$, hence 
$$\left\langle\left[
\matrix{
E_{L,T}(X_{1,1}) & \cdots & E_{L,T}(X_{1,p})\cr
\cdot & \cdots & \cdot\cr 
\cdot & \cdots & \cdot\cr 
\cdot & \cdots & \cdot\cr 
 E_{L,T}(X_{p,1}) & \cdots & E_{L,T}(X_{p,p})\cr}
\right]
\left[ 
\matrix{
x_{1} \cr 
\cdot\cr
\cdot\cr
\cdot \cr
x_{p}\cr}
\right],
\left[ 
\matrix{
x_{1} \cr 
\cdot\cr
\cdot\cr
\cdot \cr
x_{p}\cr}
\right]
\right\rangle
=\sum_{i,j=1}^{p}\langle E_{L,T}(X_{i,j})x_{j},x_{i}\rangle$$
$$=L\left(\sum_{i,j=1}^{p}(\{[l_{i,j},T^{\ast n}X_{i,j}T^{n}]p(n)^{-2}\}_{n\geq 1}\right)$$
$$
=L\left(\left\{\left\langle
\left[ 
\matrix{
X_{1,1} & \cdots & X_{1,p} \cr
\cdot & \cdots & \cdot \cr 
\cdot & \cdots & \cdot \cr 
\cdot & \cdots & \cdot \cr 
X_{p,1} & \cdots & X_{p,p}\cr}
\right] 
\left[ 
\matrix{
T^{n}x_{1} \cr
\cdot\cr
\cdot\cr
\cdot\cr 
T^{n}x_{p}\cr}
\right],
\left[ 
\matrix{
T^{n}x_{1} \cr
\cdot\cr
\cdot\cr
\cdot\cr 
T^{n}x_{p}\cr}
\right]
\right\rangle
p(n)^{-2}
\right\}_{n\geq 1}\right)\geq 0.$$
The positivity of the last term follows from the positivity of the matrix 
$[X_{i,j}]_{1\leq i,j\leq p}$ and the positivity of the state $L$.

(ii) Given any $l$ in $M_{\ast }$, we have 
$$[l,E_{L,T}(T^{\ast }XT)]=
L(\{[l,T^{\ast n+1}XT^{n+1}]p(n)^{-2}\}_{n\geq 1})$$
$$=L\left(\left\{{[l,T^{\ast n+1}XT^{n+1}]\over p(n+1)^{2}}
{p(n+1)^{2}\over p(n)^{2}}\right\}_{n\geq 1}\right).$$
Since $p$ is a gauge, we see that the sequence 
$(\left|p(n+1)/p(n)-c_{p}\right|)_{n\geq 1}$ is almost convergent to $0$. It
follows that $(\left| p(n+1)^{2}/p(n)^{2}-c_{p}^{2}\right| )_{n\geq 1}$
also almost converges to $0$. By Lemma 1 from [Ker], we get 
$$[l,E_{L,T}(T^{\ast }XT)] = c_{p}^{2}L(\{[l,T^{\ast n+1}XT^{n+1}]p(n+1)^{-2}\}_{n\geq 1})$$
$$= c_{p}^{2}L(\{[l,T^{\ast n}XT^{n}]p(n)^{-2}\}_{n\geq1})=c_{p}^{2}[l,E_{L,T}(X)]$$
and (ii) follows.

(iii) Let $A$, $B$ be two operators in $M$ commuting with $T$, we have 
$$[l,E_{L,T}(A^{\ast }XB)] =L(\{[l,T^{\ast n}A^{\ast}XBT^{n}]p(n)^{-2}\}_{n\geq 1})$$
$$=L(\{[l,A^{\ast }T^{\ast n}XT^{n}B]p(n)^{-2}\}_{n\geq 1})=[l,A^{\ast}E_{L,T}(X)B].$$
This establishes the formula.

(iv) follows immediately from (ii) and (iii).

(v) Let $l\in M_{\ast }$; using (iv), we get 
$$[l,E_{L,T}(E_{L,T}(X))] =L(\{[l,T^{\ast n}E_{L,T}(X)T^{n}]p(n)^{-2}\}_{n\geq 1})$$
$$=L\left(\left\{c_{p}^{2n}p(n)^{-2}\right\}_{n\geq 1}\right)[l,E_{L,T}(X)]=
\rho_{L}(p)[l,E_{L,T}(X)],$$
by setting $\rho _{L}(p)=L(\{c_{p}^{2n}p(n)^{-2}\}_{n\geq 1})$. From the
formula $c_{p}=\inf \{p(n)^{1/n}: n\in {\bf N}\}$ (see [Ker] Proposition 1),
we immediately deduce that $\rho _{L}(p)\in [0,1]$.

(vi) It suffices to show that the linear functional $M\ni X\mapsto l(E_{L,T}(X))$
is ultra-weakly continuous for any $l\in M_{\ast }$. Let $Z$ be in $M$, for
clarity we will denote by $l_{Z}$ the element in $M_{\ast }$ given by 
$l_{Z}(X)=\tau (ZX)$ for any $X\in M$. Given $X,Y\in M$, we have 
$$\tau (E_{L,T}(X)Y) =[l_{Y},E_{L,T}(X)]=L(\{[l_{Y},T^{\ast n}XT^{n}]p(n)^{-2}\}_{n\geq 1})$$
$$=L(\{\tau (YT^{\ast n}XT^{n})p(n)^{-2}\}_{n\geq 1})=
L(\{\tau(XT^{n}YT^{\ast n})p(n)^{-2}\}_{n\geq 1})$$
$$=[l_{X},E_{L,T^{\ast }}(Y)]=\tau (XE_{L,T^{\ast }}(Y)),$$
hence 
$$ \tau (E_{L,T}(X)Y)=\tau (XE_{L,T^{\ast }}(Y)).$$
Let $Y$ be in $M$, we deduce from the last equation that $l_Y\circ E_{L,T}$ 
is ultra-weakly continuous. Since the linear functionals $l_Y$ with $Y\in M$  are dense in $M_{\ast}$, 
it follows that $E_{L,T}$ is ultra-weakly continuous. This completes the
proof. Q.E.D.

\bigskip

{\bf II. Convergence of $\tau $-$M$-stable maps}

\medskip

Given any $X\in M$, the linear functional $l_{X}(Y):=\tau (XY)$ ($Y\in M$)
is weak*-continuous, and so $l_{X}\in M_{\ast }$. The mapping 
$$\Psi _{\tau }\colon M\to M_{\ast }, \; X\mapsto l_{X}$$
is a bounded linear quasiaffinity; the linear manifold
 $\widehat{M}_{\tau}:=\hbox{ran}\Psi $ is dense in $M_{\ast }$.

Let us consider the set 
$${B}_{\tau }(M):=\{\phi \in B_{w}(M):\phi _{\ast }(\widehat{M}_{\tau})\subset \widehat{M}_{\tau }\}$$
of $\tau $-$M$-stable weak*-continuous operators.
For any $\phi \in \widehat{B}_{\tau }(M)$, we can introduce the linear mapping 
$$\widehat{\phi }_{\tau }:=\Psi _{\tau }^{-1}\phi _{\ast }\Psi _{\tau}\colon M\to  M.$$
For any $X,Y\in M,$ we have 
$$\tau (X\phi(Y))=[l_{X},\phi(Y)]=[\phi_{\ast}(l_{X}),Y]=[l_{{\widehat{\phi }_{\tau}}(X)},Y]=
\tau (\widehat{\phi }_{\tau }(X)Y).$$
An application of the Closed Graph Theorem yields that $\widehat\phi_\tau$ is bounded.
In fact,
we see that $\widehat{\phi }_{\tau }$ is also in $\widehat{B}_{\tau }(M)$
and we have $\widehat{\left( \widehat{\phi }_{\tau }\right) }_{\tau }=\phi $.
 We say that $\phi $ is $M$-stable if it is $\tau $-$M$-stable for every $\tau \in {\cal T}(M)$.
We will consider the set 
$$\widehat{B}(M)=\bigcap_{\tau \in {\cal T}(M)}\widehat{B}_{\tau }(M)$$
of all $M$-stable operators.

\medskip

{\bf Remarks 2.1}. {\bf 1.} Denote by ${\cal M}_{p}({\bf C})$
the algebra of square matrices of order $p$, and consider the finite von
Neumann algebra $M=\oplus _{p\geq 2}{\cal M}_{p}({\bf C})$ acting on
the Hilbert space $H=\oplus _{p\geq 2}{\bf C}^{p}$ in an obvious sense.
We consider the faithful normal traces $\tau_{1}$ and $\tau_{2}$ defined
by setting 
$$\tau _{1}(\oplus_{p\geq 2}X_{p})=\sum_{p\geq 2}{1\over p^{3}}\hbox{Tr}(X_{p}),$$
$$\tau _{2}(\oplus_{p\geq 2}X_{p})=\sum_{p\geq 2}\alpha _{p}\hbox{Tr}(X_{p}),$$
where $\hbox{Tr}(\cdot)$ is the usual trace acting on ${\cal M}_{p}({\bf C})$ and $\alpha _{p}$ is given by 
$$\alpha _{p}=\left\{ 
\matrix{
{1\over2^{p}} & \hbox{ if } & p\not\in 3^{\bf N}\cr
 & & \cr
{p\over 2^p} & \hbox{ if } & p\in 3^{\bf N}.\cr}
\right.$$
 Let us consider the mapping $\phi $ defined by 
$$\phi (X_{2},X_{3},...)=(0,X_{2}',X_{3}',...)$$
where $X_{p}'\in {\cal M}_{p+1}({\bf C})$ is given in an
obvious sense by 
$$X_{p}'=\left[ 
\matrix{
X_{p} & 0 \cr
0 & 0\cr}
\right].$$
Then, we can  check that $\phi \in \widehat{B}_{\tau _{1}}(M)$ but  $\phi\not\in \widehat{B}_{\tau _{2}}(M)$.

{\bf 2.} If $M$ is a factor ($M\cap M'={\bf C}I$), then all
faithful normal traces are proportional (see [Dix, p. 249]). Consequently we
have $\widehat{B}(M)=\widehat{B}_{\tau }(M)$ for every $\tau \in {\cal T}(M)$.

{\bf 3.} Let $M$ be a finite von Neumann algebra and $A,B\in M$, then the
mapping $\phi\colon X\mapsto AXB$ belongs to $\widehat{B}(M)$.

\medskip

Let $M$ be a finite von Neumann algebra. Recall that ${\cal M}_{n}(M)$ is
also a finite von Neumann algebra with the faithful normal trace $\tau _{n}$
defined by setting 
$$\tau _{n}([A_{i,j}]_{1\leq i,j\leq n})=\sum_{k=1}^{n}\tau (A_{k,k}).$$
We begin with some useful properties of the operators $\widehat{\phi }$ when $\phi $ is a $\tau $-$M$-stable mapping.

\medskip

{\bf Proposition 2.2.}
{\it Let $M$ be a finite von Neumann
algebra, $\tau $ a faithful normal trace on $M$ and $\phi \in B(M)$. Then}

(i) {\it the mapping $\phi $ belongs to $\widehat{B}_{\tau }(M)$
 if and only if there exists $\psi \in B(M)$ such that 
$\tau (\phi (X)Y)=\tau (X\psi (Y))$ for every $X,Y\in M$; and then $\psi=\widehat\phi_\tau$};

(ii) {\it the set $\widehat{B}_{\tau }(M)$ is an algebra; the mapping $\phi\mapsto \widehat\phi_\tau$ 
is linear, involutive and
$\widehat{(\phi_1\phi_2)}_\tau=(\widehat\phi_1)_\tau(\widehat\phi_2)_\tau$};

(iii) {\it if $\phi \in \widehat{B}_{\tau }(M)$, 
then $\phi_{n}\in \widehat{B}_{\tau }({\cal M}_{n}(M))$ and we have 
$(\widehat{\phi_{n}})_{\tau }=(\widehat{\phi }_{\tau })_{n}$};

(iv) {\it if $\phi $ is $n$-positive $(n\in {\bf N})$, then $\widehat{\phi }_{\tau }$ is also $n$-positive};

(v) {\it if $\phi $ is completely positive, then 
$\widehat{\phi }_{\tau }$ is completely positive};

(vi) {\it assume that $\phi $ is $2$-positive, then
the mappings $\phi $ and $\widehat{\phi }_{\tau }$ extend
uniquely to bounded operators from $L^{2}(M,\tau )$ into itself;
moreover, we have} 

$$\left\| \phi \right\| _{B(L^{2}(M,\tau ))}=
\left\| \widehat{\phi }_{\tau}\right\| _{B(L^{2}(M,\tau ))}\leq
 \left(\left\| \phi (I)\right\|_M\right)^{1/2}\left(\left\| \widehat{\phi }_{\tau }(I)\right\| _M\right)^{1/2}.$$

\medskip

{\bf Proof.} (i) If $\phi \in \widehat{B}_{\tau }(M)$, it suffices to
set $\psi =\widehat{\phi }_{\tau }$.
 Conversely, assume that there exists 
$\psi \in B(M)$\ such that $\tau (\phi (X)Y)=\tau (X\psi (Y))$\ for every $X,Y\in M$.
 We immediately deduce that the linear functional $X\mapsto\tau(\phi (X)Y)$ is ultra-weakly continuous for each $Y\in M$. 
Since $\widehat{M}_{\tau }$ is dense in $M_{\ast }$, we see that $\phi \in B_{w}(M)$.
Moreover, we have 
$$\phi _{\ast }(l_{X})=l_{\psi (X)}$$
for any $X\in M$, thus we have $\phi_{\ast }(\widehat{M}_{\tau })\subset\widehat{M}_{\tau }$.
 This gives $\phi \in \widehat{B}_{\tau }(M)$.

(ii) This statement follows clearly from the characterization of elements of $\widehat{B}_{\tau }(M)$ given in (i).

(iii) Assume that $\phi \in \widehat{B}_{\tau }(M)$. Let $A=[A_{i,j}]_{1\leq i,j\leq n}$ 
and $B=[B_{i,j}]_{1\leq i,j\leq n}$\ be two elements in ${\cal M}_{n}(M)$, then 
$$\tau _{n}((\widehat{\phi })_{n}(A)B) =
\tau _{n}([\widehat{\phi }(A_{i,j})]_{1\leq i,j\leq n}[B_{i,j}]_{1\leq i,j\leq n})
=\sum_{k=1}^{n}\tau\left(\sum_{l=1}^{n}\widehat{\phi }(A_{k,l})B_{l,k}\right)$$
$$=\sum_{k=1}^{n}\sum_{l=1}^{n}\tau (\widehat{\phi }(A_{k,l})B_{l,k})=
\sum_{k=1}^{n}\sum_{l=1}^{n}\tau (A_{k,l}\phi (B_{l,k}))$$
$$=\tau _{n}(A\phi _{n}(B)).$$
It follows easily by (i) that $\phi _{n}\in \widehat{B}_{\tau }({\cal M}_{n}(M))$ and we have
 $\widehat{\phi _{n}}=(\widehat{\phi })_{n}$.

(iv) Assume that $\phi \in \widehat{B}_{\tau }(M)$ is positive. Let $A$ and $B$
 be two positive elements in $M$, we have 
$$\tau (\widehat{\phi }_{\tau }(A)B)=\tau (A\phi (B))=\tau (\sqrt{A}\phi (B)\sqrt{A})\geq 0.$$
Hence, we derive easily the positivity of $\widehat{\phi }_{\tau }$ from the
previous calculation. If $\phi $ is $n$-positive, the map $\phi _{n}$ is
positive, thus $(\widehat{\phi _{n}})_{\tau }$ is positive and the formula 
$(\widehat{\phi _{n}})_{\tau }=(\widehat{\phi }_{\tau })_{n}$ implies that $\widehat{\phi }_{\tau }$ is $n$-positive.

(v) It is clear from (iv) that $\widehat{\phi }$\ is completely positive if $\phi $ is completely positive.

(vi) Let $\phi \in \widehat{B}_\tau(M)$\ be $2$-positive and $Y\in M$, then the
matrix 
$$\left[
\matrix{\phi (I) & \phi (Y) \cr
\phi (Y)^{\ast } & \phi (Y^{\ast }Y)\cr}
\right]$$
is positive, a fact which implies that 
$$\phi (Y)^{\ast }\phi (Y)\leq \left\| \phi (I)\right\| \phi (Y^{\ast }Y).$$
Given a pair $(X,Y)$ of elements of $M$, we deduce from the previous
inequality that 
$$\left| \tau (\widehat{\phi }_{\tau }(X)Y)\right| =
\left| \tau (X\phi(Y))\right| \leq \sqrt{\tau (X^{\ast }X)}\sqrt{\tau (\phi (Y)^{\ast }\phi(Y))}$$
$$\leq \sqrt{\left\| \phi (I)\right\| }\sqrt{\tau (X^{\ast }X)}\sqrt{\tau(\phi (Y^{\ast }Y))}$$
$$=\sqrt{\left\| \phi (I)\right\| }\| X\| _{2}\sqrt{\tau(Y^{\ast }Y\widehat{\phi }_{\tau }(I))}$$
$$\leq \sqrt{\left\| \phi (I)\right\| }\sqrt{\left\| \widehat{\phi }_{\tau}(I)\right\| }\|X\|_{2}\|Y\|_{2}.$$
It follows that $\|\widehat\phi_\tau (X)\|_{2}\leq
 \sqrt{\left\| \phi(I)\right\| }\sqrt{\left\| \widehat{\phi }_{\tau }(I)\right\| }\|X\|_{2}$. 
Using the density of $M$ in $L^{2}(M,\tau )$, we see that
the map $\widehat{\phi }_{\tau }$\ extends uniquely to a bounded operator
from $L^{2}(M,\tau )$ into itself. We also get 
$$\left\| \widehat{\phi }\right\| _{B(L^{2}(M))}\leq
 \sqrt{\left\| \phi(I)\right\| _M} \cdot\sqrt{\left\| \widehat{\phi }_{\tau }(I)\right\| _M}.$$
Observe that the adjoint of $\phi $ in $L^{2}(M,\tau )$ is given by 
$\phi^{\ast }(X)=\widehat{\phi }_{\tau }(X^{\ast })^{\ast }$ for any $X\in M$.
The rest of the proof follows immediately. Q.E.D.

\medskip

{\bf Remark 2.3.} It follows immediately from (ii) that $\widehat{B}(M)$
is also an algebra.

\medskip

Let $\Phi =(\phi _{n})_{n\geq 1}$ be a sequence in $B(M)$ dominated by a
gauge $p$. Given a Banach limit $L$, let us consider the limit operator $E_{\Phi ,L}\in B(M)$, defined by 
$$[l,E_{\Phi ,L}(X)]=L(\{[l,\phi _{n}(X)]p(n)^{-1}\}_{n\geq 1})$$
for any $(l,X)\in M_{\ast }\times M$. Note that the previous formulas
actually define $E_{\Phi ,L}$ as an element of $B(M)$. 
We write $\gamma _{L}(p)=L(\{c_{p}^{n}p(n)^{-1}\}_n)$.

The following theorem seems to be of independent interest. It presents some
fine properties of abelian sequences included in $\widehat{B}_{\tau }(M)$
which are compatible with a gauge $p$, where $M$ is a finite von Neumann
algebra and $\tau $ is a faithful normal trace on $M$.

\medskip

{\bf Theorem 2.4.} {\it Let $M$ be a finite von Neumann
algebra, $\tau $ a faithful normal trace on $M$ and 
$\Phi =(\phi _{n})_{n\geq 1}$  a sequence in $\widehat{B}_{\tau}(M)$
 dominated by a gauge $p$ and such that 
$\widehat{\Phi }_{\tau }=((\widehat{\phi _{n}})_{\tau })_{n\geq 1}$ is also
dominated by $p$.}

(i) {\it The operator $E_{\Phi ,L}$ belongs to $\widehat{B}_{\tau }(M)$ and
we have $(\widehat{E_{\Phi ,L}})_{\tau }=E_{\widehat{\Phi }_{\tau },L}$
 for any Banach limit $L$.}

(ii) {\it Suppose $\Phi =(\phi _{n})_{n\geq 1}$ is abelian,
then the operators $E_{\Phi ,L_{1}}$ and $E_{\Phi ,L_{2}}$
commute for any pair $(L_{1},L_{2})$ of Banach limits.}

(iii) {\it Assume that $\psi \in \widehat{B}_{\tau }(M)$ and
that the sequences $\phi _{n}=\psi ^{n}$ and $\widehat{\phi _{n}}=
\widehat{\psi }^{n}$ are dominated by the gauge $p$.
Then we have} 
$$E_{\Phi ,L_{2}}\circ E_{\Phi ,L_{1}}=
E_{\Phi ,L_{1}}\circ E_{\Phi,L_{2}}=\gamma _{L_{1}}(p)E_{\Phi ,L_{2}}=\gamma _{L_{2}}(p)E_{\Phi ,L_{1}}$$
{\it for any Banach limits $L_{1}$ and $L_{2}$. In
particular, if $c_{p}^{n}/p(n)$ almost converges to a nonzero
limit, then $(p(n)^{-1}\phi _{n})_{n\geq 1}$ weakly almost
converges to an operator $\phi $ belonging to $\widehat{B}_{\tau }(M)$.}

\medskip

{\bf Proof.} (i) Given a pair $(X,Y)$ of elements of $M$, we get 
$$\tau (E_{\Phi ,L}(X)Y) =[l_{Y},E_{\Phi ,L}(X)]=L(\{[l_{Y},\phi_{n}(X)]p(n)^{-1}\}_{n\geq 1})$$
$$=L(\{\tau (\phi _{n}(X)Y)p(n)^{-1}\}_{n\geq 1})=L(\{\tau (X(\widehat{\phi_{n}})_{\tau }(Y))p(n)^{-1}\}_{n\geq 1})$$
$$=L(\{[l_{X},(\widehat{\phi _{n}})_{\tau }(Y)]p(n)^{-1}\}_{n\geq1})
=[l_{X},E_{\hat{\Phi }_{\tau },L}(Y)]=\tau (XE_{\hat{\Phi }_{\tau},L}(Y))$$
for every Banach limit $L$.
Now Proposition 2.2.(i) implies the statement.

(ii) Assume that the sequence $(\phi _{n})_{n\geq 1}$ is abelian. We have 
$$ [l_{(\widehat{\phi _{m}})_{\tau }(Y)},\phi _{n}(X)]=
\tau (\phi_{n}(X)(\widehat{\phi _{m}})_{\tau }(Y))=
\tau (\phi _{m}(X)(\widehat{\phi_{n}})_{\tau }(Y))=[l_{\phi _{m}(X)},(\widehat{\phi _{n}})_{\tau }(Y)]$$
for any pair $(m,n)$ of positive integers and any pair $(X,Y)\in M^{2}$.
By taking  $L_{1}$-limit  with respect to  $n$  we get 
$$[l_{E_{\Phi ,L_{1}}(X)},(\widehat{\phi _{m}})_{\tau }(Y)] =
\tau(E_{\Phi ,L_{1}}(X)(\widehat{\phi _{m}})_{\tau }(Y))=
[l_{(\widehat{\phi _{m}})_{\tau }(Y)},E_{\Phi ,L_{1}}(X)]$$
$$=
[l_{\phi _{m}(X)},E_{\hat{\Phi }_{\tau },L_{1}}(Y)]=
\tau (\phi_{m}(X)E_{\hat{\Phi }_{\tau },L_{1}}(Y))=
[l_{E_{\hat{\Phi }_{\tau},L_{1}}(Y)},\phi _{m}(X)].$$
Now, taking $L_{2}$-limit with respect to $m$ and using (i), we obtain that 
$$\tau (E_{\Phi ,L_{2}}\circ E_{\Phi ,L_{1}}(X)Y)
=\tau (E_{\Phi,L_{1}}(X)(\widehat{E_{\Phi ,L_{2}}})_{\tau }(Y))
=\tau(E_{\Phi ,L_{1}}(X)E_{\hat{\Phi }_{\tau },L_{2}}(Y))$$
$$=[l_{E_{\Phi ,L_{1}}(X)},E_{\hat{\Phi }_{\tau },L_{2}}(Y)] 
=[l_{E_{\hat{\Phi },L_{1}}(Y)},E_{\Phi ,L_{2}}(X)]
=\tau (E_{\hat{\Phi }_{\tau },L_{1}}(Y)E_{\Phi ,L_{2}}(X))$$
$$=\tau ((\widehat{E_{\Phi ,L_{1}}})_{\tau }(Y)E_{\Phi ,L_{2}}(X))
=\tau(E_{\Phi ,L_{1}}\circ E_{\Phi ,L_{2}}(X)Y).$$
We thus have $E_{\Phi ,L_{2}}\circ E_{\Phi ,L_{1}}=E_{\Phi ,L_{1}}\circ E_{\Phi ,L_{2}}$.

(iii) Let $\psi \in B_{w}(M)$ and $(L_{1},L_{2})$ be a pair of Banach
limits. For any $X,Y\in M$, we have 
$$[l_{\phi _{n}(X)},(\widehat{\phi _{m}})_{\tau }(Y)]p(n)^{-1}p(m)^{-1}
=\tau (\phi _{n}(X)(\widehat{\phi _{m}})_{\tau}(Y))p(n)^{-1}p(m)^{-1}$$
$$=\tau (\psi ^{m+n}(X)Y)p(n)^{-1}p(m)^{-1}={p(m+n)\over p(m)p(n)}[l_{Y},\psi ^{m+n}(X)]p(m+n)^{-1}.$$
The sequence $p(m+n)/p(m)=p(m+1)/p(m)...p(m+n)/p(m+n-1)$ is strongly almost
convergent to $c_{p}^{n}$, when $m$ goes to infinity.
 By taking $L_{1}$-limit  with respect
to $m$, we get 
$$[l_{E_{\hat{\Phi }_{\tau },L_{1}}(Y)},\phi _{n}(X)]p(n)^{-1}
=\tau (\phi _{n}(X)E_{\hat{\Phi }_{\tau },L_{1}}(Y))p(n)^{-1}$$
$$=[l_{\phi _{n}(X)},E_{\hat{\Phi }_{\tau},L_{1}}(Y)]p(n)^{-1}=
c_{p}^{n}p(n)^{-1}[l_{Y},E_{\Phi ,L_{1}}(X)].$$
After taking  $L_{2}$-limit  with respect to $n$, it follows that 
$$\tau(E_{\Phi ,L_{1}}\circ E_{\Phi ,L_{2}}(X)Y)=
[l_{E_{\hat{\Phi }_{\tau },L_{1}}(Y)},E_{\Phi ,L_{2}}(X)]$$
$$=
\gamma _{L_{2}}(p)[l_{Y},E_{\Phi,L_{1}}(X)]
=\gamma _{L_{2}}(p)\tau (E_{\Phi ,L_{1}}(X)Y),$$
whence $E_{\phi,L_1}\circ E_{\phi,L_2}=\gamma_{L_2}(p)E_{\phi,L_1}$.
We deduce from (ii) that $E_{\phi,L_{1}}$ and $E_{\phi,L_{2}}$ commute.
Interchanging the role of $L_{1}$ and $L_{2}$ we 
conclude that
 $$E_{\Phi ,L_{2}}\circ E_{\Phi ,L_{1}}=
E_{\Phi ,L_{1}}\circ E_{\Phi ,L_{2}}=\gamma _{L_{1}}(p)E_{\Phi ,L_{2}}=\gamma _{L_{2}}(p)E_{\Phi,L_{1}}.$$

When the sequence $c_{p}^{n}/p(n)$ almost converges to a nonzero number, then we deduce that the
limit of $(p(n)^{-1}\phi _{n})_{n\geq 1}$ is independent of $L$. Applying
Lorentz's result, we obtain that $(p(n)^{-1}\phi _{n})_{n\geq 1}$ is weakly
almost convergent to an operator in $\widehat{B}_{\tau }(M)$. The proof is
now complete.
Q.E.D.

\medskip

{\bf Remarks 2.5.} {\bf 1.} If $c_{p}^{n}p(n)^{-1}$ is almost
convergent to a nonzero limit, then we see that the Cesaro means 
$(n+1)^{-1}(p(0)^{-1}\phi_{0}+...+p(n)^{-1}\phi _{n})$  weakly converge to an operator in 
$B_{w}(M)$ (which is obviously ultra-weakly continuous).

{\bf 2.} The assumption that $\Phi =(\phi _{n})_{n\geq 1}$\ is a sequence
in $\widehat{B}_{\tau }(M)$\ dominated by a gauge $p$\ does not imply that 
$\widehat{\Phi }_{\tau }=((\widehat{\phi _{n}})_{\tau })_{n\geq 1}$ is also
dominated by $p$. Denote by ${\cal M}_{p}({\bf C})$ the algebra of
square matrices of order $p$, and consider the finite von Neumann algebra 
$M=\oplus _{p\geq 2}{\cal M}_{p}({\bf C})$ acting on the Hilbert space 
$H=\oplus _{p\geq 2}{\bf C}^{p}$ in an obvious sense. We consider the
faithful normal trace $\tau $ defined by 
$$\tau (\oplus _{p\geq 2}{X}_{p})=
\sum_{p\geq 2}{1\over p^{3}}\hbox{Tr}({X}_{p}),$$
where $\hbox{Tr}(\cdot)$ is the usual trace  on ${\cal M}_{p}({\bf C})$. 
Fix a unit vector $e_{p}$ in ${\bf C}^{p}$ and write $P_{p}$ for the orthogonal
projection onto ${\bf C}^{p}\ominus {\bf C}e_{p}$. 
For any $n\geq 2$,
set $\phi _{n}(\oplus _{p\geq 2}X_{p})=
\langle X_{n}e_{n}, e_{n}\rangle P_{n}$. 
We can easily see that $\Phi =(\phi _{n})_{n\geq 1}$\ is a\ sequence in 
$\widehat{B}_{\tau }(M)$ such that $\left\| \phi _{n}\right\| =1$ for every $n$.
 Hence $\Phi $ is dominated by the constant gauge $p$ equal to $1$, but
$\widehat{\Phi }_{\tau }=((\widehat{\phi _{n}})_{\tau })_{n\geq 1}$ is not dominated
by $p$ , actually $\|(\widehat{\phi _{n}})_{\tau }\| =n-1$.

\medskip

Let ${\cal F}$ be an abelian set included in $\widehat{B}_{\tau }(M)$.
We consider the (abelian) semigroup  ${\cal E(F)}$ induced by $\cal F$, that is
$${\cal E(F)}=\{\phi _{_{1}}\circ\dots\circ \phi _{n}: \phi _{_{1}},...,\phi_{n}\in {\cal F}, n\in {\bf N}\}.$$
We define a partial ordering on ${\cal E(F)}$ by setting 
$F\le F'$ if there exists $F''$ in $\cal E(F)$ such that $F'=F''F$.
(It is clear that $\cal E(F)$ is a directed set with this partial ordering, and it can be considered as a
net (generalized sequence) indexed by itself.)

\medskip

{\bf Proposition 2.6. }
{\it  Let $M$ be a finite von Neumann algebra, 
$\tau $ a faithful normal trace on $M$, and 
let ${\cal F}$ be an abelian set of $m$-positive projections
belonging to $\widehat{B}_{\tau }(M)$ ($m\geq 2$).
 Assume
that ${\cal E(F)}$ and $\widehat{{\cal E(F)}}$ are
bounded in $B(M)$.
 Then the net ${\cal E(F)}$ weakly
converges to an $m$-positive projection $E\in \widehat{B}_{\tau}(M)$.}

\medskip

{\bf Proof.}  Let us introduce the classical Hilbert space $L^{2}(M,\tau) $ 
equipped with the inner product 
$\langle X,Y\rangle:=\tau(Y^*X)\; (X,Y\in M)$.
 Since every $\phi\in {\cal E(F)}$ is $2$-positive, Proposition 2.2.(vi) shows that $\phi $
extends uniquely to a bounded projection (still denoted by $\phi $) from $L^{2}(M,\tau )$ into itself.
 Using again Proposition 2.2.(vi), we see that
the set ${\cal E(F)}$ is bounded in $B(L^{2}(M,\tau ))$, thus it is
weakly relatively compact.

Choose two cofinal subnets $(E_{i})_{i\in {\cal I}}$ and $(E_{j})_{j\in {\cal J}}$ in ${\cal E(F)}$,
 which converge respectively to $F$ and $G$ in the weak operator
topology of $B(L^{2}(M,\tau ))$. 
Fix $i\in {\cal I}$ and consider the set ${\cal J}_{i}=\{j\in {\cal J}: j\geq i\}$.
Then we have 
$$\langle E_{i}(X), E_{j}^{\ast }(Y)\rangle 
=\tau (E_{i}(X)(\widehat{E_{j}})_{\tau }(Y^{\ast }))=\tau (E_{j}\circ E_{i}(X)Y^{\ast })$$
$$=\tau (E_{j}(X)Y^{\ast })=\langle E_{j}(X), Y\rangle$$
for any $j\in {\cal J}_{i}$ and any pair $(X,Y)\in M^{2}$.
 Thus, taking limit with respect to the set ${\cal J}_{i}$, we obtain 
$$\langle E_{i}(X), G^{\ast }(Y)\rangle=\langle G(X), Y\rangle.$$
Now, taking limit with respect to the directed  set $\cal I$, we get 
$$\langle G\circ F(X), Y\rangle=\langle F(X), G^{\ast }(Y)\rangle=\langle G(X), Y\rangle,$$
whence $G\circ F=G$ follows.
 Interchanging the role of $F$ and $G$, we see that 
$F\circ G=F$.
Since $F$ and $G$ are limit points of elements belonging to the commutative
set ${\cal E(F)}$, they commute.
Hence $F=F\circ G=G\circ F=G$, in particular $F=F\circ G=F\circ F$.
We deduce that ${\cal E(F)}$ is weakly convergent in $B(L^{2}(M,\tau ))$ to a projection 
$E$.

Now, we want to show that $E(M)\subset M$ and $E|M\in B(M)$.
To this order let $X\in M$ be arbitrary, and let us consider the linear functional
$\varphi(Y):=\langle E(X),Y^*\rangle \; (Y\in M)$.
Choosing a cofinal subnet $\{E_i\}_{i\in{\cal I}}$ in $\cal E(F)$,
we have $\langle E(X),Y^*\rangle=\lim_i\langle E_i(X),Y^*\rangle=\lim_i\tau(E_i(X)Y)$.
Thus
$$|\varphi(Y)|=\lim_i|\tau(E_i(X)Y)|\le\liminf_i\|E_i(X)\| \|Y\|_1\le C\|X\| \|Y\|_1,$$
where $C=\sup\{\|F\|: F\in{\cal E(F)}\}<\infty$ and
$\|Y\|_1:=\tau(|Y|)=\|l_Y\|$ (see [Dix] Section I.6.10).
We deduce that there exists unique $\tilde X\in M$ such that
$$\langle E(X),Y^*\rangle=\varphi(Y)=(\varphi\circ\Psi_\tau^{-1})(l_Y)=[l_Y,\tilde X]=
\tau(\tilde XY)=\langle\tilde X,Y^*\rangle$$
holds for every $Y\in M$.
Hence $E(X)=\tilde X\in M$ and $\|E(X)\|=\|\tilde X\|\le C\|X\|$.

It is clear that
$$[l_Y,E(X)]=\tau(E(X)Y)=\langle E(X),Y^*\rangle=\lim_i\langle E_i(X),Y^*\rangle=\lim_i[l_Y,E_i(X)]$$
is true for every $X,Y\in M$.
Since $\widehat M_\tau$ is dense in $M_*$, and $\cal E(F)$ is bounded, it follows that $\cal E(F)$ weakly converges to $E$.

We can prove in the same manner that $\widehat{\cal E(F)}$ converges weakly to an operator $F\in B(M)$.
Taking into account that $\tau(E_i(X)Y)=\tau(X\widehat{(E_i)}_\tau(Y))$, we obtain by passing to the limit
that 
$$\tau(E(X)Y)=\tau(XF(Y))$$
holds for every $X,Y\in M$. It follows by Proposition 2.2.(i) that $E$ belongs to $B_w(M)$.

It remains to prove that $E$ is $m$-positive.
It is clear that every operator in $\cal E(F)$ is $m$-positive.
Let $[X_{k,l}]_m\in{\cal M}_m(M)$ be a positive operator.
Given any vector $x=x_1\oplus\cdots\oplus x_m\in H^{(m)}$, we have
$$\left\langle\left[E(X_{k,l})\right]_m x,x\right\rangle=
\sum_{k,l=1}^m \left[ l_{x_l,x_k},E(X_{k,l})\right]$$
$$=\lim_i[l_{x_l,x_k},E_i(X_{k,l})]=\lim_i\left\langle\left[E_i(X_{k,l})\right]_m x,x\right\rangle\ge0,$$
and so $E$ is an $m$-positive projection.
Q.E.D.

\medskip

Let $M$ be a finite von Neumann algebra and let $T$ be an operator in $M$
dominated by the regular gauge $p$.
We know by
 Theorem 2.4.(iii)
 that the sequence of $M$-stable mappings 
$(\phi _{T,n})_{n\geq 0}$
defined by $\phi _{T,n}(X)=p(n)^{-2}T^{\ast n}XT^{n}$ is weakly almost
convergent, we will denote its limit by $E_{T}$.
Now Proposition 1.2  shows that $E_{T}$
is a completely positive projection. Notice also that $E_{T}$ is an
ultra-weakly continuous $M$-stable operator with $\widehat E_T=E_{T^*}$ and that $\|E_T\|\le1$.
Let ${\cal S}$ be an abelian subset of $M$,
 which is dominated by the  regular gauge $p$.
 We consider the abelian semigroup $\cal E(S)$ induced by the (abelian) set
$\{E_T: T\in{\cal S}\}$.

\medskip

{\bf Corollary 2.7. } {\it Let $M$ be a finite von Neumann
algebra and let $\cal S$ be an abelian subset of $M$.
Assume that $\cal S$ is dominated by a regular gauge $p$.
Then the net ${\cal E(S)}$ weakly converges to a completely
positive $M$-stable projection $E$ satisfying the
following properties}:

(i) {\it $E(T^{\ast }XT)=c_p^2 E(X)$ for any $T\in \cal S$ and
$X\in M$};

(ii) {\it $E(A^{\ast }XB)=A^{\ast }E(X)B$ for any pair $(A,B)\in ({\cal S}')^2$ and $X\in M$.}

\medskip

{\bf Proof.}  We apply Proposition 2.6 to the net ${\cal E(S)}$. 
We
deduce that ${\cal E(S)}$\ converges weakly to $E\in \widehat{B}_{\tau}(M)$, for any $\tau\in{\cal T}(M)$. 
Since properties (i) and (ii) are
true for the operators $E_R \; (R\in\cal S)$ by Proposition 1.2, we see that the same properties hold for $E$.
Q.E.D.

\medskip

{\bf Proposition 2.8.}
{\it Let $M$ be a finite von Neuman
algebra, $\tau $ a faithful normal trace on $M$, and let $\cal S$ be an
abelian subset of $M$ which is dominated by a
regular gauge $p$.
 Let us assume that the limit projection $E$
 of the net ${\cal E(S)}$ is such that $E(I)$ is injective.
Then there exists an abelian set ${\cal S}_{1}$ of
unitaries belonging to $M$ such that for any $T\in \cal S$
 there exists $U_{T}\in {\cal S}_1$ satisfying $\sqrt{E(I)}T=c_pU_{T}\sqrt{E(I)}$. 
Moreover, if $F$ denotes
the limit of the net ${\cal E(S}_{1}{\cal )}$, then we have the
following properties}:

(i) {\it $E(\sqrt{E(I)}X\sqrt{E(I)})=\sqrt{E(I)}F(X)\sqrt{E(I)}$
 for any $X\in M$};

(ii) {\it $\sqrt{E(I)}\widehat{E}_{\tau }(X)\sqrt{E(I)}=
\widehat{F}_{\tau }(\sqrt{E(I)}X\sqrt{E(I)})$ for any $X\in M$.}

\medskip

{\bf Proof.} (i) First of all, observe that the equation $T^{\ast}E(I)T=c_{p}^{2}E(I)$ and the injectivity of $E(I)$ 
imply that $\sqrt{E(I)}T$
is also injective.
Taking the polar decomposition of $\sqrt{E(I)}T$, we see
that there exists a unique isometry $U_{T}$ such that $\sqrt{E(I)}T=c_{p}U_{T}\sqrt{E(I)}$. 
Since $U_T\in M$ and $M$ is finite, it follows that $U_T$ is unitary.
The previous intertwining relations readily imply that the set
${\cal S}_1:=\{U_T: T\in {\cal S}\}$ is abelian.

Given $T\in {\cal S}$ and $X\in M$, we have 
$$p(n)^{-2}T^{\ast n}\sqrt{E(I)}X\sqrt{E(I)}T^{n}=
{c_{p}^{2n}\over p(n)^{2}}\sqrt{E(I)}U_{T}^{\ast n}XU_{T}^{n}\sqrt{E(I)}$$
for every positive integer $n$. Taking a Banach limit we get the relation
$$E_T(\sqrt{E(I)}X\sqrt{E(I)})=\sqrt{E(I)}E_{U_T}(X)\sqrt{E(I)}.$$
Now taking limits in the nets $\cal E(S)$ and ${\cal E(S}_1)$ we get (i).

(ii) Since $E$ and $F$ are $\tau $-$M$-stable, we can now get (ii) by the
following computation. For any $(X,Y)\in M^{2}$ we have
$$\tau (X\sqrt{E(I)}\widehat{E}_{\tau }(Y)\sqrt{E(I)})
=\tau (\sqrt{E(I)}X\sqrt{E(I)}\widehat{E}_{\tau }(Y))
=\tau (E(\sqrt{E(I)}X\sqrt{E(I)})Y)$$
$$=\tau (\sqrt{E(I)}F(X)\sqrt{E(I)}Y)
=\tau (X\widehat{F}_{\tau }(\sqrt{E(I)}Y\sqrt{E(I)})).$$
Hence, we have 
$\sqrt{E(I)}\widehat{E}_{\tau }(Y)\sqrt{E(I)}=\widehat{F}_{\tau }(\sqrt{E(I)}Y\sqrt{E(I)})$ for any $Y\in M$. 
This completes the proof.
Q.E.D.

\bigskip

{\bf III. Similarity}

\medskip

We say that an operator $T$ is {\it asymptotically controlled} by a gauge $p$ if $T $ is compatible with $p$
 and satisfies the  condition that $q'(\{\|T^nx\|^2/p(n)^2\}_n)>0$, for every nonzero vector $x\in H$,
where
$$q'(\xi):=\sup\left\{\liminf_k  {1\over m}\sum_{i=1}^m \xi(n_i+k): m\in{\bf N}, n_1,\dots,n_m\in{\bf N}\right\}$$
for any bounded real sequence $\xi$
(see [Ker] for the role of this functional in the study of Banach limits).

For any  real sequence $\xi\in\ell^\infty({\bf N}^n) \; (n\in{\bf N}, n>1)$,
let $Q_n\xi:=\eta\in\ell^\infty({\bf N}^{n-1})$, where
$\eta(j_1,\dots,j_{n-1}):=q'(\xi_{j_1,\dots,j_{n-1}})$ with
$\xi_{j_1,\dots,j_{n-1}}(j):=\xi(j_1,\dots,j_{n-1},j)$.
Let $\widetilde Q_n:=Q_1\circ \dots \circ Q_{n-1}\circ Q_n$, where $Q_1:=q'$.

A set $\cal F$ of operators, acting on the Hilbert space $H$, is called asymptotically controlled by a gauge $p$, if
every operator in $\cal F$ is compatible with $p$, and if for every nonzero vector $x\in H$
there exists $\rho(x)>0$ such that
$$\widetilde Q_n\left(\left\{{1\over p(j_1)^2\cdots p(j_n)^2}\left\langle T_n^{*j_n}\cdots T_1^{*j_1}T_1^{j_1}\cdots T_n^{j_n}x,
x\right\rangle\right\}_{j_1,\dots,j_n=1}^\infty\right)\ge\rho(x)$$
is true for every $n\in{\bf N}$ and $T_1,\dots,T_n\in{\cal F}$.

\medskip

{\bf Remark 3.1.} Let $T$ be an operator compatible with a gauge $p$.
Assume that $T$ satisfies  
$$\inf \{\left\| T^{n}x\right\| /p(n): n\in {\bf N}\}>0$$
for any nonzero $x$ in $H$; then $T$ is asymptotically controlled by $p$.
 In particular, power bounded operators of class $C_{1\cdot}$ (in the Sz.-Nagy--Foias
terminology) are exactly operators which are asymptotically controlled by
constant gauges.

\medskip

{\bf Theorem 3.2.}
{\it Let $\cal S$ be an abelian set of operators
which is contained in a finite von Neumann algebra. 
Assume that $\cal S$ is asymptotically controlled by a regular gauge $p$. 
Then, there
exists an invertible operator $A$ in $M$ such that 
$r(T)^{-1}ATA^{-1}$ is a unitary operator for any $T\in \cal S$.}

\medskip

{\bf Proof. } Let $\tau $ be a faithful normal trace acting on $M$.
Let $E$ be the completely positive limit projection provided by Corollary 2.7.
 Since $\cal S$ is asymptotically controlled by the gauge $p$, we can infer by a short computation that,
 given any nonzero vector $x\in H$,
$$[l_{x,x},E_{T_1}\circ E_{T_2}\circ\cdots\circ E_{T_n}(I)]\ge \rho(x)$$
is true for  every choice of $T_1,\dots,T_n\in{\cal S}, n\in{\bf N}$ with a $\rho(x)>0$,
whence 
$$\langle E(I)x,x\rangle=[l_{x,x},E(I)]\ge\rho(x)>0.$$
Thus, the positive operator $E(I)$ is injective.
Let us consider the associated set ${\cal S}_1$ and the corresponding limit operator $F$ occurring in Proposition 2.8.

Set $X=E(I)$, $Y=\widehat{E}_{\tau }(I)$ and consider the positive operator $R=\sqrt{X}Y\sqrt{X}=\widehat{F}_{\tau }(X)$.
Note that $R$ commutes with $U_T\; (T\in{\cal S})$.
Let $P$ be a projection
associated with the spectral decomposition of $R$ (which still commutes with $U_T$).
 By the Cauchy--Schwarz
Inequality, we get 
$$\tau (PRP)=\tau (P\sqrt{X}Y\sqrt{X}P)\leq \sqrt{\tau (PXP)}\sqrt{\tau (P\sqrt{X}Y^{2}\sqrt{X}P)}.\leqno(1)$$
Applying  the properties of $E$ and $F$ described in Corollary 2.7 and
Proposition 2.8, we infer that 
$$\tau (PXP) =\tau (XP)=\tau (XF(P))=\tau (\sqrt{X}F(P)\sqrt{X})=\tau (E(\sqrt{X}P\sqrt{X}))$$
$$=\tau (\sqrt{X}P\sqrt{X}\widehat{E}_{\tau }(I))=\tau (\sqrt{X}P\sqrt{X}Y)=\tau (PR)=\tau (PRP).$$
Now, note that the operator $C=\sqrt{X}P\sqrt{X}Y$ commutes with $T^{\ast }$, 
because we have
$$T^*C=T^*\sqrt X P\sqrt X Y= c_p\sqrt X U_T^* P\sqrt X Y=c_p\sqrt X PU_T^*\sqrt X Y$$
$$=(1/c_p)\sqrt X PU_T^*\sqrt X TYT^*=\sqrt X PU_T^*U_T\sqrt X YT^*=\sqrt XP\sqrt X YT^*=CT^*.$$
Using again Corollary 2.7 and Proposition 2.8 we obtain 
$$\tau (P\sqrt{X}Y^{2}\sqrt{X}P) =\tau (Y\sqrt{X}P\sqrt{X}Y)=\tau (YC)=\tau (\widehat{E}_{\tau }(I)C)
=\tau (\widehat{E}_{\tau }(C))$$
$$=\tau (\widehat{E}_{\tau }(C)E(I))=\tau (\widehat{E}_{\tau }(I)CE(I))=\tau(Y\sqrt{X}P\sqrt{X}YX)$$
$$=\tau (P\sqrt{X}YXY\sqrt{X})=\tau (PR^{2})=\tau (PR^{2}P).$$
Substituting these results into (1), we get 
$$\tau (PRP)\leq \sqrt{\tau (PRP)}\sqrt{\tau (PR^{2}P)},$$
whence
$$\tau (PRP)\leq \tau (PR^{2}P).\leqno(2)$$
Let $K$ be a compact set contained in the interval $(0,1)$, and denote by $P$ the
spectral projection associated to $K$ by the functional calculus of $R$.
We
thus have 
$$  PRP\ge (PRP)^{2}=PR^{2}P\ge0.\leqno(3)$$
Combining (2) with (3) yields 
$$\tau (PRP-PR^{2}P)=0.$$
The operator $PRP-PR^{2}P$ is positive, so it is necessarily equal to $0$.
Therefore $Q=PRP$ is an orthogonal projection.
 But $K$ is compact and
contained in $(0,1)$, thus there exists $\rho \in (0,1)$ such that 
$Q\leq\rho I$. Consequently, we have $Q=0$. 
It follows that 
$$\sigma (R)\cap (0,1)=\emptyset.$$

The last step is devoted to show that $0\not\in \sigma (R)$. 
Let us denote by $P$
the spectral projection associated to $0$. 
We have $RP=0$, thus 
$$0=\tau (RP)=\tau (\widehat{F}_{\tau }(X)P)=\tau (XF(P))=\tau (XF(I)P)=\tau(XP)=\tau (PXP).$$
It follows that $PXP=0$. Since $X$ is injective, we deduce that $P=0$. 
Finally,
we see that $R$ is invertible (actually, $\sigma (R)\subset [1,\infty)$),
 therefore $X$ is also invertible. 
By Proposition 2.8 we know that $\sqrt X T=r(T)U_T\sqrt X \; (c_p=r(T)$, see [Ker]).
It follows that, for any $T\in\cal S$, $\; U_T=r(T)^{-1}ATA^{-1}$ is unitary, where $A=\sqrt X$ is invertible.
Q.E.D.

\bigskip

{\bf Acknowledgments}. The author wishes to thank the organizers of the
``Memorial Conference for B\'ela Sz\H okefalvi-Nagy", the organizers of 
``Journ\'{e}es d'Analyse Fonctionnelle de Lens'' and the University of Besan\c con where these results were presented.
The author wishes to express his gratitude to L\'aszl\'o K\'erchy for some
stimulating questions.

\bigskip

{\bf References}

\bigskip

\noindent

\item{[Bea]} B. Beauzamy, {\it Introduction to Operator Theory and Invariant
Subspaces}, North-Hol-land, Amsterdam (1988).

\medskip
\noindent

\item{[BCP]} S. Brown, B. Chevreau and C. Pearcy, Contractions with rich
spectrum have invariant subspaces, {\it J. Operator Theory}, 
{\bf 1} (1979), 123--136.

\medskip
\noindent

\item{[CaFa]} G. Cassier and T. Fack, Contractions in von Neumann algebras,
{\it J. Funct. Anal.}, {\bf 135} (1996),  297--338.

\medskip
\noindent

\item{[CaFa2]} G. Cassier and T. Fack, On power-bounded operators in finite
von Neumann algebras, {\it J. Funct. Anal.}, {\bf 141} (1996),  133--158.

\medskip
\noindent

\item{[Dix]} J. Dixmier, {\it Les alg\`{e}bres d'op\'{e}rateurs dans l'espace
Hilbertien (Alg\`{e}bre de von Neumann)}, Gauthier-Villards, Paris (1969).

\medskip
\noindent

\item{[Dix2]} J. Dixmier,  Les moyennes invariantes dans les semi-groupes et
leurs applications, {\it Acta Sci. Math.} ({\it Szeged}), {\bf 12} (1950),  213--227.

\medskip
\noindent

\item{[Lor]} G. G. Lorentz,  A contribution to the theory of divergent
sequences, {\it  Acta. Math.}, {\bf 80} (1948),  167--190.

\medskip
\noindent

\item{[Ker]} L. K\'erchy, Operators with regular norm sequences, {\it Acta Sci.
Math.} ({\it Szeged}), {\bf 63} (1997),  571--605.

\medskip
\noindent

\item{[Ker2]} L. K\'erchy, Criteria of regularity for norm sequences,
{\it Integral Equations Operator Theory}, {\bf 34} (1999),  458-477.

\medskip
\noindent

\item{[Ker3]} L. K\'erchy, Representations with regular norm-behavior of
discrete abelian semigroups, {\it Acta Sci. Math.}  ({\it Szeged}), {\bf 65} (1999),  702--726.

\medskip
\noindent

\item{[Ker4]} L. K\'erchy, Hyperinvariant subspaces of operators with
non-vanishing orbits, {\it Proc. Amer. Math. Soc.}, {\bf 127} (1999),  1363--1370.

\medskip
\noindent

\item{[Ker5]} L. K\'erchy, Unbounded representations of discrete abelian
semigroups, {\it Progress in Nonlinear Differential Equations and Their
Applications}, {\bf 42} (2000),  141--150.

\medskip
\noindent

\item{[KeM\"{u}]} L. K\'erchy and V. M\"{u}ller, Criteria of regularity for
norm sequences. II, {\it Acta Sci. Math.}  ({\it Szeged}), {\bf 65} (1999),  131--138.

\medskip
\noindent

\item{[Sak]} S. Sakai, {\it $C^{\ast }$ algebras and $W^{\ast }$
algebras}, Ergebnisse der Mathematik und ihrer Grenzgebiete, {\bf 60} (1971),
Springer-Verlag, Berlin--New York.

\medskip
\noindent

\item{[Nag]} B. Sz.-Nagy, On uniformly bounded linear transformations in
Hilbert space, {\it Acta Sci. Math.} ({\it Szeged}), {\bf 11} (1947),  152--157.

\bigskip

Institut Girard Desargues

UPRES-A 5028 Math\'ematiques

Universit\'e Claude Bernard Lyon I

69622 Villeurbanne Cedex, France

{\it e-mail}: cassier@desargues.univ-lyon1.fr

\bye